# Geometric mean, splines and de Boor algorithm in geodesic spaces

Esfandiar Nava-Yazdani

**Abstract**

We extend the concepts of de Casteljau and de Boor algorithms as well as splines to geodesic spaces and present some applications in geometric modeling. The concept of weighted geometric mean provides another approach to splines. We compare the corresponding Bézier curves and show that for Riemannian manifolds, their endpoint tangents coincide.

## 1 Introduction

A geodesic space is a metric space where any two points can be joined by a shortest geodesic realizing the distance between those points. If there is a unique shortest geodesic between any two points (unique geodesic spaces), affine combination of points is well defined and in turn, the natural and general setting for de Casteljau and de Boor algorithms is provided. This construction shares many interesting properties with the special case of manifold-valued data considered in [6], [17], [29] and [25] (as well as [32], [26] and [31] for subdivision schemes). Moreover, geodesic spaces give rise to a wider range of applications. A survey on particular geodesic spaces including Alexandrov and Busemann spaces can be found in [18]. A recent related work, considering subdivision schemes in metric spaces and link to barycenters in Hadamard spaces, is [10]. Furthermore, an extrinsic approach to spline curves in embedded manifolds by minimizing energy, and several significant examples can be found in [14] and [15]. Data refinement in nonlinear geometries is a substantial problem in geometric processing and has a wide range of applications. Furthermore, many properties and applications of Bernstein polynomials immediately extend to analogous constructions in geodesic spaces. For example, an $\mathbb{R}^m$-valued function $f$ can be approximated by the polynomial $\sum_{i=0}^{n} f(i/n) B_i^n$ where $B_i^n$ denotes the $i$-th Bernstein polynomial of degree $n$. Similar results for manifold-valued functions are desirable. For Bernstein polynomials on spheres and sphere-like surfaces we refer to [3]. The intrinsic approach in the present work can be used to construct

Bézier curves, and more generally, splines. Usually, to ensure well definedness of our approach, restriction of control points to a small enough neighbourhood is sufficient. In Riemannian manifolds, the size of this neighbourhood (for the algorithmic construction of splines to be well defined) is determined by the injectivity radius. In the present work we focus on geometric modeling.

The concept of geometric mean is important in many applications such as medical imaging and numerical linear algebra (cf. [24] and [12]). Moreover, it plays a significant role in studying some essential geometric and analytic properties. Besides the classic work [20], we would like to refer to [1], [19], [9] and [8] for some recent progresses and applications. For a comprehensive treatment of spherical averages and splines, we refer to [5]. Choosing the weights as the Bernstein polynomials, the corresponding wighted geometric mean curves, provide an approach to splines. We present some basic results concerning the existence and uniqueness of the weighted geometric mean, lower bound estimates of the cost functional, and briefly compare the corresponding curves with Bézier curves. In particular, it turns out, that their tangents at endpoints coincide.

This work is organized as follows. In the next section, we present some preliminaries on geodesic spaces. The second section is devoted to construction of Bézier curves via the de Casteljau algorithm in unique geodesic spaces, rational Bézier curves, the de Boor algorithm, splines and some examples. The last section is devoted to the concept of weighted geometric mean, the resulting curves for the case that the weights in the cost functional are the Bernstein polynomials and an application of the Karcher equation.

## 2 Preliminaries

We recall a few definitions and refer to [18] for details. Let $(M, d)$ be a metric space. The length $L$ of a curve $c \in C^0([0, 1], M)$ is

$$L(c) := \sup \left\{ \sum_{i=0}^{n-1} d(c(t_i), c(t_{i+1})) \, : \, 0 = t_0, \cdots, t_n = 1, \, n \in \mathbb{N} \right\}.$$

$c$ is called geodesic iff there exists $\epsilon > 0$ with

$$L(c|_{[s,t]}) = d(c(s), c(t)) \quad \text{whenever } |s - t| < \epsilon.$$

Its image is called geodesic arc or segment and will be denoted by $[c(0), c(1)]$. A geodesic $c : [0, 1] \to M$ is called shortest geodesic iff $L(c) = d(c(0), c(1))$.

**Definition 2.1.** *We call $(M, d)$ a (unique) geodesic space[1] iff for any two points in $M$ there exists a (unique) shortest geodesic arc joining them. If $(M, d)$ is a unique geodesic space, than we define*

$$\Phi \, : \, M \times M \times [0, 1] \ni (x, y, t) \mapsto \Phi_t(x, y) \in M$$

*as the unique point of the shortest geodesic arc between $x$ and $y$ with $d(x, \Phi_t(x, y)) = td(x, y)$ and refer to $\Phi$ as the affine map of $M$.*

---
[1] Some authors use the terminology geodesic length space.

In this setting, the notation of betweenness also makes sense:
$$y \in [x, z] \Leftrightarrow d(x, y) + d(y, z) = d(x, z).$$

For example, any neighbourhood of a complete Riemannian manifold within the injectivity radius (in particular, any Hadamard-Cartan manifold) is a unique geodesic space. Euclidean trees, and more generally, Bruhat-Tits buildings provide examples of nonmanifold unique geodesic spaces. Denoting the Euclidean inner product by $\langle .,. \rangle$ for any $x, y$ in an open hemisphere of $S^2$ we have

$$\Phi_t(x, y) = \frac{\sin((1-t)\varphi)}{\sin \varphi} x + \frac{\sin(t\varphi)}{\sin \varphi} y \text{ with } \varphi := \arccos(\langle x, y \rangle).$$

Note that in general for geodesic spaces, restriction to a small neighbourhood, does not result in a unique geodesic space. For instance, denoting the Euclidean norm by $|.|$, every neighbourhood in $\mathbb{R}^2$ endowed with the Manhattan (Taxicab) metric defined by

$$d(x, y) = |x_1 - y_1| + |x_2 - y_2|,$$

is geodesic but not unique. Nevertheless, fixing a representative of identified geodesics, an affine map can be simply defined. For example, let $k \in \mathbb{R}$, $c = (x + y)/2$, $g$ the line passing through the point $c$ with slope $k$ and $x^*$ resp. $y^*$ the nearest point of $g$ to $x$ resp. $y$. Consider the geodesic $[x, x*] \cup [x*, y*] \cup [y*, y]$ with $[ , ]$ being Euclidean. Obviously, the corresponding affine map reads

$$\Phi_t^k(x, y) = \begin{cases} (1 - \frac{L_3 t}{L_1})x + \frac{L_3 t}{L_1} x^* & \text{for } 0 \leq t < L_1/L_3, \\ \frac{L_2 - L_3 t}{L_2 - L_1} x^* + \frac{L_3 t - L_1}{L_2 - L_1} y^* & \text{for } L_1/L_3 \leq t \leq L_2/L_3, \\ \frac{L_3 - L_3 t}{L_3 - L_2} y^* + \frac{L_3 t - L_2}{L_3 - L_2} y & \text{for } L_2/L_3 < t \leq 1, \end{cases}$$

where $L_1 = |x - x^*|$, $L_2 = L_1 + |x^* - y^*|$ and $L_3 = L_2 + |y^* - y|$.

Throughout this work, unless otherwise stated explicitly, $(M, d)$ will denote a unique geodesic space and $\Phi$ its affine map. Furthermore, we refer to $p_0, \cdots, p_n \in M$ as control points and set $I := [0, 1]$.

## 3 Bézier and spline curves

Now, we present the definition of Bézier curves via the generalized de Casteljau algorithm in unique geodesic spaces.

### 3.1 De Casteljau algorithm and Bézier curves

**Definition 3.1.** *Fix $t \in I$. Let $p_i^0(t) := p_i$, $i = 0, \cdots, n$ and set*

$$r = 1, \cdots, n,$$
$$i = 0, \ldots, n - r,$$
$$p_i^r(t) := \Phi_t(p_i^{r-1}(t), p_{i+1}^{r-1}(t)).$$

*We call $I \ni t \mapsto p(t) := p_0^n(t)$ the Bézier curve with control points $p_0, \cdots, p_n$.*

Obviously this curve is invariant under affine parameter transformations, it lies in the convex hull
$$\{\Phi_I([p_i, p_{i+1}], [p_j, p_{j+1}]) : i, j = 0, \cdots, n-1\}$$
of the control points and has the end point property $p(0) = p_0$, $p(1) = p_n$. Moreover, replacing $t$ by $\frac{t-t_i}{t_{i+r}-t_i}$ for $0 = t_0 < t_1 < \cdots < t_n = 1$, we get the Aitken-Neville interpolation of $p(t_i) = p_i$ in unique geodesic spaces.

Applications of de Casteljau algorithm to some Riemannian manifolds including Lie groups and, more generally, symmetric spaces can be found in [6], [17], [29] and [25]. Next we consider a metric real tree.

**Example 3.1.** *Denote $[x, y, z] = 0$ iff $x, y$ and $z$ lie in the same geodesic. Let $c \in M$. The Paris metric $d_p$ with respect to $c$ can be defined as*
$$d_p(x, y) = \begin{cases} d(x, y) & \text{if } [x, y, c] = 0, \\ d(x, c) + d(y, c) & \text{else.} \end{cases}$$

*Obviously, we have $\Phi^{paris}(x, y) = \Phi(x, y)$ if $[x, y, c] = 0$, and*
$$\Phi_t^{paris}(x, y) = \begin{cases} \Phi_{\frac{Lt}{L_1}}(x, c) & \text{for } 0 \leq t \leq L_1/L, \\ \Phi_{\frac{Lt-L_1}{L-L_1}}(y, c) & \text{for } L_1/L \leq t \leq 1 \end{cases}$$

*otherwise. Here $L := d(x, c) + d(y, c)$ and $L_1 := d(x, c)$. Utilizing the above exlicit formulae for the affine map $\Phi^{paris}$, corresponding Bézier curves can be easily computed. We would just like to mention that Bézier curves for certain combinations of the Paris (with more complicated trees), Manhattan and Euclidean metric serve in many industrial applications like path generation and planning (for instance in urbanistic, automated cartography and mobile robots construction, cf. [2], [23],[11] and [22]), and postpone a consideration of those topics to another opportunity.*

Imposing a further condition on $M$, Bézier curves defined in 3.1 enjoy the following subdivision property.

**Theorem 3.1.** *Suppose that*
$$\Phi_s(\Phi_\tau(x, y), \Phi_\tau(y, z)) = \Phi_\tau(\Phi_s(x, y), \Phi_s(y, z)), \quad s, \tau \in I, \ x, y, z \in M. \tag{1}$$

*Then, for any partition $0 \leq s_1 \leq \cdots \leq s_k \leq 1$, the Bézier curve $p$ can be split into $k+1$ Bézier curves on $[0, s_1], \cdots, [s_k, 1]$.*

*Proof.* Let $k = 1$ and $s \in ]0, 1[$. Let $x_i := p_0^i(s)$, $i = 0, \cdots, n$. We show by induction on $n$ that the left segment $I \ni t \mapsto p^{left}(t) := p_0^n(st)$ of the Bézier curve $p_0^n$ and the Bézier curve $x := x_0^n$, determined by the control points $x_i$, coincide. The identity
$$\Phi_t(z_0, \Phi_s(z_0, z_1)) = \Phi_{st}(z_0, z_1), \quad s, t \in I, \ z_0, z_1 \in M$$

is crucial for the proof. Fix $t \in I$. Let $n = 1$. Then $x_0^1(t) = \Phi_t(x_0, x_1) = \Phi_t(p_0, \Phi_s(p_0, p_1)) = \Phi_{st}(p_0, p_1) = p_0^1(st) = p^{left}(t)$. Now, suppose that $n > 1$ and the assertion is true for $n-1$ control points. Hence $x_0^{n-1}(t) = p_0^{n-1}(st)$. Moreover, the point $x_1^{n-1}(t)$ is determined by the control points $p_0^1(s) = \Phi_s(p_0, p_1), \cdots, p_0^n(s) = \Phi_s(p_0^{n-1}(s), p_1^{n-1}(s))$. Hence

$$x_1^{n-1}(t) = \Phi_t(p_0^{n-1}(s), \Phi_s(p_0^{n-1}(s), p_1^{n-1}(s))) = \Phi_{st}(p_0^{n-1}(s), p_1^{n-1}(s)).$$

Applying (1) (with $\tau = st$), we get $x_1^{n-1}(t) = \phi_s(p_0^{n-1}(st), p_1^{n-1}(st))$. Therefore

$$\begin{aligned} x_0^n(t) &= \Phi_t(x_0^{n-1}(t), x_1^{n-1}(t)) = \Phi_t(p_0^{n-1}(st), \Phi_s(p_0^{n-1}(st), p_1^{n-1}(st))) \\ &= \Phi_{st}(p_0^{n-1}(st), p_1^{n-1}(st)) = p_0^n(st). \end{aligned}$$

Hence $p^{left} = x$. Reversing the order of the control points and replacing $s$ by $1-s$ and $t$ by $1-t$, the same argument shows that the right segment $I \ni t \mapsto p^{right}(t) := p_0^n(s+t-st)$ of the Bézier curve $p_0^n$ and the Bézier curve $y := y_0^n$, determined by the control points $y_i = p_i^{n-i}(s)$, are equal. Hence, we arrive at $p([0,1]) = x([0,s]) \cup y([s,1])$, $x([0,s]) \cap y([s,1]) = \{p(s)\}$. We iterate on $k$ to complete the proof. □

Condition (1) seems to be natural for the compatibility between subdivision property and de Casteljau iteration and could be of interest for further investigation. The following theorem gives further properties of the de Casteljau algorithm. Proofs are slight modifications of the Riemannian case presented in [25] and we outline them for the reader's convenience.

**Theorem 3.2.** *Consider control points $P := (p_0, \ldots, p_n)$ and $Q := (q_0, \ldots, q_n)$ in a compact neighbourhood $U \subset M$.*
*a) Transformation invariance: Suppose a Lie group $H$ acts on $M$ by*

$$H \times M \ni (h, x) \mapsto hx \in M$$

*leaving $U$ invariant, i.e., $HU \subset U$. Denote $B(p_0, \cdots, p_n) := p(I)$. If the action is segment-equivariant, i.e., for every $h \in G$*

$$h[x, y] = [hx, hy] \text{ for all } x, y \in U.$$

*Then*

$$hB(p_0, \ldots, p_n) = B(hp_0, \ldots, hp_n) \text{ for all } h \in H.$$

*b) Local control:*

$$\sup_{t \in I} d(p(t), q(t)) \leq C d_\infty(P, Q) := |d(p_0, q_0), \cdots, d(p_n, q_n)|_\infty$$

*where $C$ denotes a positive constant depending only on $n$ and $U$.*

*Proof.* a) The segments $[hp_i, hp_{i+1}]$ and $h[p_i, p_{i+1}]$ have the same endpoints:

$$\Phi_0(hp_i, hp_{i+1}) = hp_i = h\Phi_0(p_i, p_{i+1}),$$
$$\Phi_1(hp_i, hp_{i+1}) = hp_{i+1} = h\Phi_1(p_i, p_{i+1}).$$

b) Fix $t \in I$. There is a positive constant $K$ determined by $U$ and Lipschitz constant of $\Phi$ on $U$ such that

$$d(p_i^r(t), q_i^r(t)) = d(\Phi_t(p_i^{r-1}(t), p_{i+1}^{r-1}(t)), \Phi_t(q_i^{r-1}(t), q_{i+1}^{r-1}(t)))$$
$$\leq K|(d(p_i^{r-1}(t), q_i^{r-1}(t)), d(p_{i+1}^{r-1}(t)), q_{i+1}^{r-1}(t)))|_\infty.$$

Iteration yields

$$d(p(t), q(t)) = d(p_0^n(t), q_0^n(t)) \leq K^n d_\infty(P, Q)$$

which immediately implies the desired inequality. □

## 3.2 Rational Bézier curves

In a geodesic space $(M, d)$, the rational Bézier curve $p_0^n$ with control points $p_0, \cdots, p_n$ and positive weights $w_0, \cdots, w_n$ can be produced by applying the weighted version of the de Casteljau algorithm as follows. Fix $t \in I$. Let $(p_i^0(t), w_i^0(t)) := (p_i, w_i)$, $i = 0, \cdots, n$ and set

$$r = 1, \ldots, n,$$
$$i = 0, \ldots, n - r,$$
$$w_i^r(t) = (1 - t)w_i^{r-1}(t) + tw_{i+1}^{r-1}(t),$$
$$t_i^r = t\frac{w_{i+1}^{r-1}(t)}{w_i^r(t)}, \ p_i^r(t) := \Phi_{t_i^r}(p_i^{r-1}(t), p_{i+1}^{r-1}(t)).$$

The following figures show the effect of weights. We remark that in general, for a surface of revolution, the geodesic differential equation reduces to a first order one and can efficiently be solved using e.g. ode45 of MATLAB.

**Example 3.2.** *In order to reflect constraints caused by the presence of some objects (without changing the degree of the Bézier curve), we may choose weights as functions of distances from control points to the objects. In this example we treat the small disc $\mathcal{B}$ in $M$ as an attracting object. Here each weight is simply chosen as inverse of the distance between corresponding control point and $\mathcal{B}$, i.e., $w_i = 1/d(p_i, \mathcal{B})$. Similarly, avoiding objects can be treated. Of course, due to the convex hull property of Bézier curves, the gained flexibility is restricted. For a treatment of obstacles via a variational approach we refer to [15] and [14].*

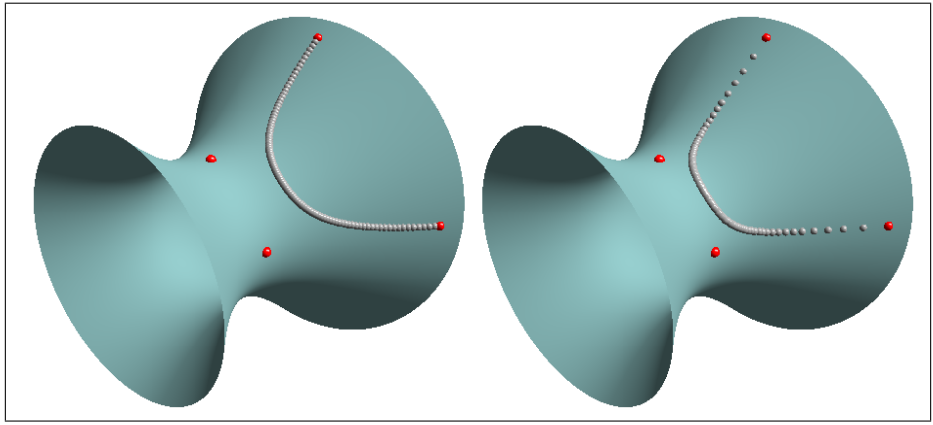

Figure 1: Bézier curve in a one-sheeted hyperboloid: left) cubic , right) rational with weights $1, 5, 5, 1$.

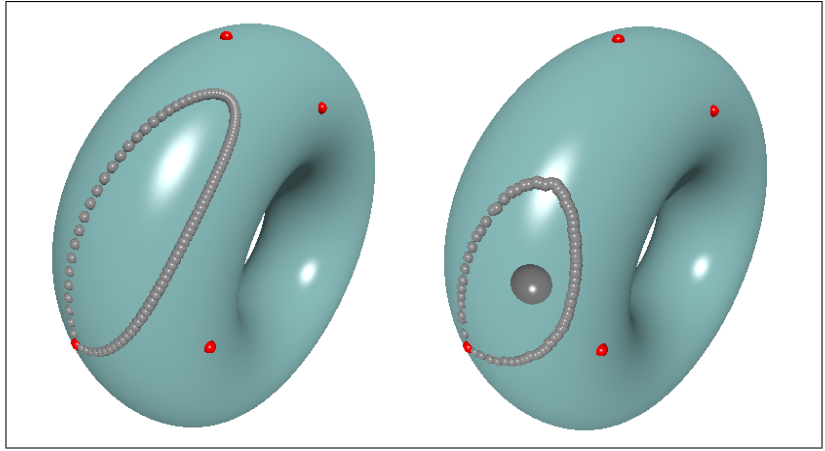

Figure 2: Bézier curve in a torus: left) cubic , right) rational with weights equal to inverse of distances between the small disc and control points.

## 3.3 De Boor algorithm and splines

For a knot vector
$$\tau : \tau_0 \leq \tau_1 \leq \cdots \leq \tau_{m+n+1}$$
with $\tau_k < \tau_{k+n}$, $\tau_m < \tau_{m+1}$, $\tau_n < \tau_{n+1}$ and control points $p_0, \cdots, p_n$ in a geodesic space $M$ with metric $d$, we define the de Boor algorithm as follows. Fix $t$ in a knot interval $[\tau_l, \tau_{l+1}[$ of the parameter interval $[\tau_m, \tau_{n+1}]$, define $p_i^0(t) := p_i$, $i = l - m, \cdots, l$ and set

$$r = 1, \ldots, m,$$
$$i = r - m + l, \ldots, l,$$
$$t_{i,\tau}^r = \frac{t - \tau_i}{\tau_{i+m-r} - \tau_i}, \; p_i^r(t) := \Phi_{t_{i,\tau}^r}(p_{i-1}^{r-1}(t), p_i^{r-1}(t)).$$

Then the spline curve $p$ of degree $\leq m$ evaluated at $t$ is obtained as the final value $p_l^m$. If $t = \tau_l$ with multiplicity $\mu$, then the the scheme simplifies slightly and $p(t) = p_{l-\mu}^{m-\mu}$. Note that due to continuity, closed knot intervals are permitted. For closed splines, control points and knots should be extended periodically.

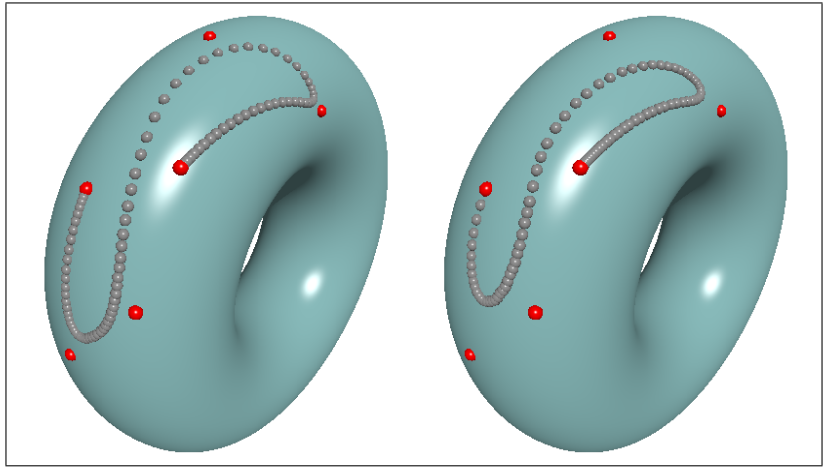

Figure 3: Cubic spline curve in a torus: left) uniform, right) with double knots.

In general, it is not easy to obtain geodesics which are required for the construction of the splines using the above algorithm. For Riemannian symmetric spaces and in particular, Lie groups, explicit formulas for the Riemannian exponential map can be used to reduce the task to calculations with matrix exponential and prinicipal logarithm accordingly. Next example illustrates these considerations.

**Example 3.3.** *Poses of a rigid body can be visualized as a curve in the Euclidean motion group*

$$E_3 = \left\{ \begin{pmatrix} 1 & 0 \\ b & R \end{pmatrix} : R \in SO(3), b \in \mathbb{R}^3 \right\}$$

*for which the convexity radius is*[2] $\pi/2$ *and (cf. [32])* $\Phi_t(x,y) = \exp(tx \log((x^{-1}y))$.

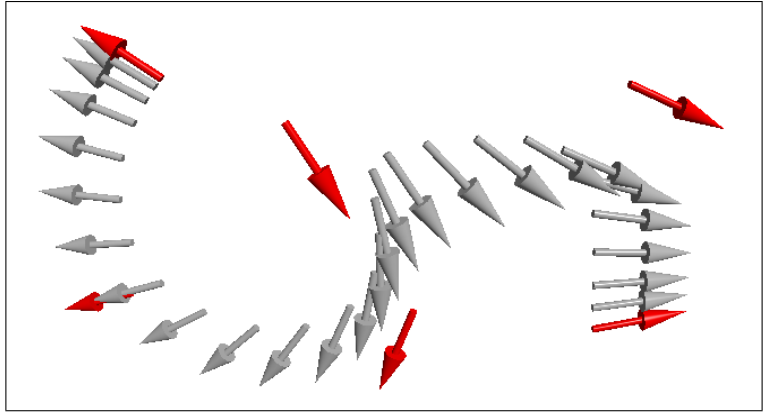

Figure 4: Cubic uniform spline curve in Euclidean motion group.

---

[2]cf. [21] and [4].

# 4 Geometric mean and centroid curves

## 4.1 Geometric mean

We recall that a function $\rho : M \to \mathbb{R}$ is said to be (strictly) convex iff $\rho \circ c$ is (strictly) convex for any nonconstant geodesic $c$. Now, let $b_i^n$ be some nonnegative real weights satisfying (w.l.o.g.) $\sum_{i=0}^n b_i^n = 1$ and consider the following optimization task

$$\sum_{i=0}^n b_i^n d^2(.,p_i) \to \min. \qquad (2)$$

on $M$. If there is a unique minimizer, then it is called (weighted) center of mass, geometric mean[3] or centroid. Throughout the remainder of this work, $B_i^n$ denotes the $i$-th Bernstein polynomial of degree $n$. In the following, we present a general existence and uniqueness result for the geometric mean and some estimates of the corresponding cost functional for the case that weights are given by the Bernstein polynomials. We say that the control points are in general position unless they are contained in a common geodesic.

**Theorem 4.1.** *Fix $t \in I$. Let $E^r := \sum_{i=0}^r B_i^r(t) d^2(., p_i^{n-r}(t))$ with $r = 1, \ldots, n$.*
*a) For any $x \in M$, we have*

$$E^n(x) \geq (t(1-t))^n (\sum_{i=0}^n \binom{n}{i} d(x,p_i))^2 \geq (t(1-t))^n (\sum_{i=0}^{n-1} \binom{n-1}{i} d(p_i,p_{i+1}))^2.$$

*For $t \in ]0,1[$, inequalities are strict if and only if the control points are in general position.*
*b) Suppose that $d^2(.,x)$ is strictly convex for all $x \in M$. Then the weighted geometric mean defined by (2) exists and is unique. Denoting $q^r := argMin(E^r)$, for any $n > 1$, the sequence $\{E^r(q^r)\}_{r=1,\ldots,n}$ is increasing. It is strictly increasing, provided $t \in ]0,1[$.*

*Proof.* a) Note that for any $i = 0, \ldots, n-1$, we have

$$E_i^1(x) := (1-t)d^2(x,p_i) + td^2(x,p_{i+1}) \geq t(1-t)(d(x,p_i) + d(x,p_{i+1}))^2$$
$$\geq t(1-t)d^2(p_i,p_{i+1}) = E_i^1(p_i^1(t))$$

where inequalities become equality if and only if $x = p_i^1(t)$. Iteration completes the proof. Note that for $n > 1$, the control points are in general position if and only if for at least one $i \in \{1, \ldots, n-1\}$, the triangle $p_{i-1}p_ip_{i+1}$ is not degenerated, implying strictness of the asserted inequalities.

b) As a positive linear combination of proper strictly convex functionals $d^2(.,p_i)$, the cost functional defining the weighted geometric mean in (2) is proper and strictly convex. Therefore it attains its unique minimum. Moreover, for any $r = 1, \ldots, n-1$, denoting

---
[3]also known as Fréchet or Karcher mean

$e_i := d^2(x, p_i^{n-r-1}(t))$ with $i = 0, \ldots, r+1$, we have

$$E^{r+1}(x) = \sum_{i=0}^{r+1} B_i^{r+1}(t)e_i = B_0^{r+1}(t)e_0 + \sum_{i=1}^{r} B_i^{r+1}(t)e_i + B_{r+1}^{r+1}(t)e_{r+1}$$

$$= B_0^{r+1}(t)e_0 + \sum_{i=1}^{r}((1-t)B_i^r(t) + tB_{i-1}^r(t))e_i + B_{r+1}^{r+1}(t)e_{r+1}$$

$$= \sum_{i=0}^{r} B_i^r(t)((1-t)e_i + te_{i+1})$$

$$\geq \sum_{i=0}^{r} B_i^r(t)d^2(x, p_i^{n-r}(t)) = E^r(x).$$

Hence $E^{r+1}(q^{r+1}) \geq E^r(q^{r+1}) \geq E^r(q^r)$. The asserted strictness of the monotonicity for $t \in ]0, 1[$, follows from the strictness of the convexity. □

This result can be used to determine the geodesic $p_0^1$ by minimizing $E^1$. Conversely, the preceding theorem can be applied to determine the geometric mean of two points, if geodesics are known. More generally, the given lower bounds can be utilized to approximate the minimizer of $E^n$ with a Casteljau-like algorithm.

**Example 4.1.** *Let $M$ be the space of symmetric positive definite $2 \times 2$-matrices of determinant one and $x, y \in M$. Then, for any $s \in I$, we have (cf. [32])*

$$\Phi_s(x, y) = x(x^{-1}y)^s.$$

*Due to the preceding theorem the geometric mean $z$ of $x$ and $y$ is given by*

$$z = \Phi_{1/2}(x, y) = x(x^{-1}y)^{1/2}.$$

*An application of Cayley-Hamilton's theorem (see [24]) implies $z = \frac{x+y}{\sqrt{det(x+y)}}$.*

## 4.2 Centroid curves

**Definition 4.1.** *Suppose that for any $x \in M$, the function $d^2(., x)$ is strictly convex. We call $I \ni t \mapsto q(t) := argMin_{x \in M}(\sum_{i=0}^{n} B_i^n(t)d^2(., p_i))$ the centroid curve of the control points $p_0, \ldots, p_n$.*

Obviously, the centroid curve is continuous and enjoys the endpoint property $q(0) = p_0$ and $q(1) = p_1$. Due to the endpoint property, several centroid curves can be pieced together into a spline. While the evaluation of $p$ relies on computation of geodesics, for $q$ we require (just) geodesic distances, but as a drawback, the solution of an optimization task (resp. Karcher equation in Riemannian case) which is in general highly nonlinear and complex.

We have seen that $p = q$, if $n = 1$. Next, we expose a brief comparison corresponding to this fact, when the weighted geometric mean is replaced by the median $argMin(F)$ where $F := (1-t)d(.,p_i) + td(.,p_1)$ with some fixed $t \in I$. Thus, suppose that $d$ is convex. Fix $x \in M$. Let $\mathcal{B}_i$ denote the ball with radius $r_i = d(x,p_i)$ around $p_i$ and $l := d(p_0,p_1)$. If $r_i > l$ for $i = 1$ or $0$, then we have $F(x) > l \geq F(q)$ for any $q \in [p_0,p_1]$. Hence, $x$ cannot be a minimizer. Therefore, we may and do assume that $r_i \leq l$. In view of the triangle inequality, there is a point $q \in [p_0,p_1]$ in the intersection of $\mathcal{B}_0$ and $\mathcal{B}_1$. Thus, $d(x,p_i) \geq d(q,p_i)$ implying $F(x) \geq F(q) = F(\Phi_s(p_0,p_1))$ for some $s \in I$ with $q = \Phi_s(p_0,p_1)$, and we arrive at $F(x) \geq ((1-t)s + t(1-s))l$. Therefore, $F$ has a minimizer in $[p_0,p_1]$. Indeed, $d(x,.)$ is linearly increasing on $[p_0, p_0^1(1/2)]$ and linearly decreasing on $[p_0, p_0^1(1/2)]$ for any $x$ and we have (breakdown at $t = 1/2$)

$$argMin(F) = \begin{cases} p_0 & \text{for } 0 \leq t < 1/2, \\ p_0^1(1/2) & \text{for } t = 1/2, \\ p_1 & \text{for } 1/2 < t \leq 1. \end{cases}$$

Next, we assume that $(M,d)$ is a smooth complete Riemannian manifold without conjugate points and consider the characterization of the geometric mean as the unique solution of the corresponding Karcher equation given below. Note that, within the injectivity radius of the Riemannian exponential map exp, in terms of local representatives, the affine map of $(M,d)$ is given by

$$\Phi_t(x,y) = \exp_x(t \log_x(y))$$

where $\log_x$ denotes the local inverse of the exponential map at $x$. We recall two important well-known special cases when for any control points the weighted geometric mean is well-defined. First, if the sectional curvatures of $M$ are bounded above by $k > 0$ and $diam(M) < \pi/(2\sqrt{k})$, second, Cartan-Hadamard manifolds like the space of positive definite symmetric matrices. Many progresses concern the latter (for which there are also recent extensions to the infinite-dimensional setting presented in [19]). For further applications as well as the proof (based on Jacobi field estimates) of the corresponding Karcher's result we refer to [1], [9], [8] and the pioneer work [20]. Some related results on the geometric median and an extension of the Weiszfeld algorithm for its computation can be found in [12].

In the following we suppose that $M$ is contained in an open ball of radius less than $\frac{1}{2}\min(r_{inj}, \frac{\pi}{2\sqrt{k}})$ of a Riemannian manifold with injectivity radius $r_{inj}$ and sectional curvatures bounded above by $k$ where $\frac{1}{\sqrt{k}} := \infty$ if $k \leq 0$. Then, due to [20], for any $y \in M$, the squared distance function $d^2(y,.)$ is strictly convex and $\sum_{i=1}^n b_i^n d^2(.,p_i)$ has a unique minimizer $x \in M$ determined by the Karcher equation

$$\sum_{i=0}^n b_i^n \log_x p_i = 0.$$

We present for $n = 1$ an alternate direct proof of the fact that $p(t) = p_0^1(t)$ is the unique minimizer of $(1-t)d^2(.,p_0) + td^2(.,p_1)$ utilizing the Karcher equation

$$(1-t)\log_x p_0 + t\log_x p_1 = 0.$$

Denoting $v := \log_{p_0} p_1$ (the initial velocity of the geodesic from $p_0$ to $p_1$), we have $\log_{p(t)} p_0 = -tw(t)$ and $\log_{p(t)} p_1 = (1-t)w(t)$ where $w(t)$ denotes the parallel transport of $v$ along $p$ at $t$. Therefore

$$(1-t)\log_{p(t)} p_0 + t\log_{p(t)} p_1 = -(1-t)tw(t) + t(1-t)w(t) = 0.$$

Hence, the point $p(t)$ is critical. With $d^2(., p_0)$ and $d^2(., p_1)$, the weighted sum $(1-t)d^2(., p_0) + td^2(., p_1)$ is also strictly convex, implying that the critical point $p(t)$ is its unique minimizer.

Next, we consider smoothness of the centroid curve and compute its tangents at endpoints.

**Theorem 4.2.** *The curve $q$ is smooth and its velocity $\dot{q}(t) := \frac{d}{d\tau}\big|_{\tau=t} q(\tau)$ satisfies*

$$\dot{q}(0) = n\log_{p_0} p_1, \ \dot{q}(1) = -n\log_{p_n} p_{n-1}.$$

*Proof.* Smoothness of the centroid curve $q$ is an application of the implicit function theorem to the Karcher equation. To prove the formula for the velocities, let $x_0 \in M$. Taking the tangent maps of both sides of the identity $\exp(\log x_0) = x_0$ and evaluating at $x_0$, yields $Id_{T_{x_0}M} + (d\log x_0)\big|_{x_0} = 0_{x_0}$, i.e., $(d\log x_0)\big|_{x_0} = -Id_{T_{x_0}M}$ where $T_{x_0}M$ denotes the tangent space at $x_0$ and $0_{x_0}$ its zero, and we used the natural identification $T_v T_{x_0}M \simeq T_{x_0}M$ for any $v \in T_{x_0}M$. Now, by the Karcher equation

$$0 = \sum_{i=0}^{n} \left( \dot{B}_i^n(0)) \log_{q(0)} p_i + B_i^n(0)(d\log_q p_i)\big|_{q=q(0)} \dot{q}(0) \right)$$
$$= n\log_{p_0} p_1 - \dot{q}(0).$$

Reversing the order of the control points and replacing $t$ by $1-t$ completes the proof. □

The preceding result can be used for $C^1$ interpolation and path fitting with centroid curves. We remark that due to [30], $\dot{p}(0) = n\log_{p_0} p_1$ and $\dot{p}(1) = -n\log_{p_n} p_{n-1}$. It follows that $p$ and $q$ have the same tangents at their endpoints. Obviously, in Euclidean spaces, the curves $p$ and $q$ coincide. We have seen that in general geodesic spaces, they are equal, if $n = 1$. In the next example we compare $p$ and $q$ for the spherical case and, as an application of the Karcher equation, show that already for three control points, even the images of those curves (although close) are different.

**Example 4.2.** *Let $p_0, p_1, p_2$ be vertices of an equilateral spherical triangle with the common length of the sides $\alpha \in ]0, \pi/2]$. A straightforward calculation shows*

$$p(1/2) = \frac{1}{4\cos(\theta/2)\cos(\alpha/2)}(p_0 + 2p_1 + p_2)$$

*where*

$$\cos\theta = \frac{1 + 3\cos\alpha}{2 + 2\cos\alpha}.$$

*Now, $p(1/2) \in q(I)$ holds, if and only if there exists $s \in I$ such that $p(1/2)$ satisfies the Karcher equation, i.e.,*

$$L(s) := \sum_{i=0}^{2} B_i^2(s) \frac{\psi_i}{\sin \psi_i}(p_i - p(1/2) \cos \psi_i) = 0$$

*where*

$$\cos \psi_1 = \langle p_1, p(1/2) \rangle = \frac{2 + 2\cos\alpha}{4\cos(\theta/2)\cos(\alpha/2)},$$

$$\cos \psi_0 = \cos \psi_2 = \langle p_0, p(1/2) \rangle = \frac{1 + 3\cos\alpha}{4\cos(\theta/2)\cos(\alpha/2)} = \cos\theta \cos\psi_1.$$

*Note that $\cos\psi_0 < \cos\psi_1$, hence $z := \frac{\psi_1 \sin\psi_0}{\psi_0 \sin\psi_1} < 1$. Denoting the cross product in $\mathbb{R}^3$ by $\times$, we have*

$$\langle L(s), p_2 \times p_0 \rangle = \left( \sum_{i=0}^{2} B_i^2(s) \frac{\psi_i \cos\psi_i}{2\cos(\theta/2)\cos(\alpha/2)\sin\psi_i} - B_1^2(s)\frac{\psi_1}{\sin\psi_1} \right) \langle p_1, p_0 \times p_2 \rangle$$

$$= \left( (1-s)^2 + s^2 - 2s(1-s)z \right) \frac{\cos\theta \langle p_1, p_0 \times p_2 \rangle \psi_0}{(1+\cos\theta)\sin\psi_0}$$

$$= \underbrace{\left( (1+z)(s-\tfrac{1}{2})^2 + \frac{1-z}{4} \right)}_{\geq \frac{1-z}{4} > 0} \frac{\cos\theta \langle p_1, p_0 \times p_2 \rangle \psi_0}{(1+\cos\theta)\sin\psi_0}.$$

*In particular, $L$ has no zero in $I$. Therefore $p(1/2) \notin q(I)$.*

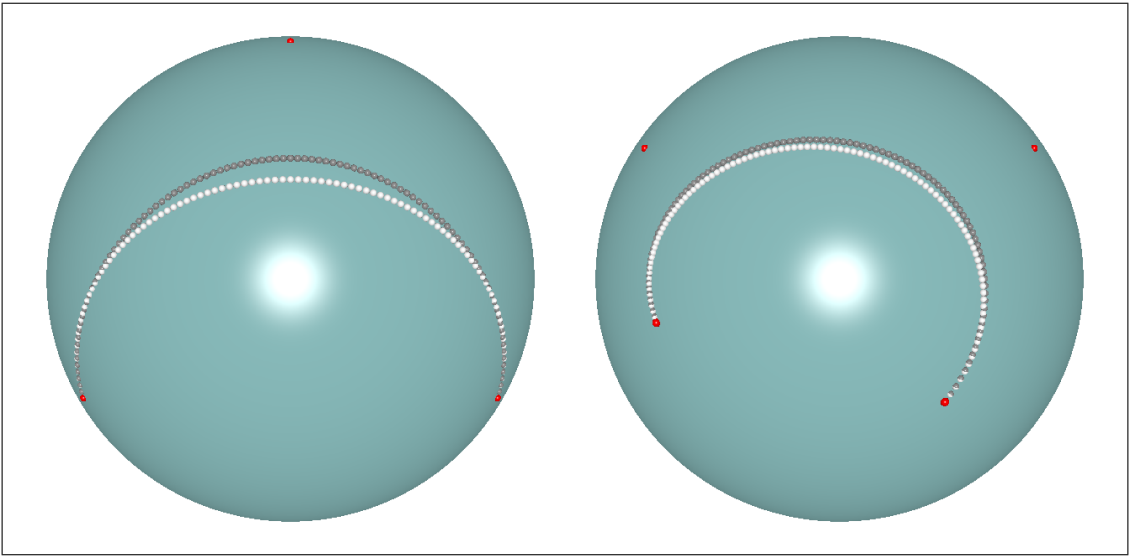

Figure 5: Spherical Bézier (gray) and centroid curves: left) quadratic with equidistant control points, right) cubic.

## 5  Conclusion

In this paper we extended the setting for de Casteljau algorithm as well as de Boor algorithm to unique geodesic spaces and presented some examples. In turn, Bézier and more generally spline curves in these spaces can be constructed via iteration. Of course, the main issue remains the fact that one needs to compute geodesics. We expect besides Riemannian manifolds, applications of our results in other geodesic spaces including trees and Bruhat-Tits buildings. Furthermore, we expect similar results concerning subdivision schemes as well as the bivariate case to produce nets and spline surfaces in geodesic spaces. Moreover, weighted extended B-splines (WEB-splines, [13] and [16]) can be combined with our approach to provide efficient approximations with high accuracy for solutions of boundary value problems on manifolds. Moreover, considering the Bernstein polynomials as weights for the cost function defining the geometric mean, and piecing them together, we get a different approach to interpolation tasks and construction of splines. We proved that the Bézier interpolation of the corresponding curves is $C^1$. Further properties of the latter curves and their application in interpolation would be desirable future tasks.